\documentclass[prb,endfloats]{revtex4}
\usepackage{graphicx}
\newcommand{\bp}{{\bf p}}
\newcommand{\bq}{{\bf q}}

\newcommand{\bac}{{\bf a}}
\newcommand{\bff}{{\bf F}}

\newcommand{\A}{{\bf A}}
\newcommand{\ct}{{\cal T}}

\newcommand{\bx}{{\bf x}}
\newcommand{\bv}{{\bf v}}

\newcommand{\ep}{h}

\newcommand{\e}{{\rm e}}

\newcommand{\nn}{\nonumber}
\newcommand{\la}{\label} 

\newcommand{\be}{\begin{equation}}
\newcommand{\ee}{\end{equation}}
\newcommand{\ba}{\begin{eqnarray}}
\newcommand{\ea}{\end{eqnarray}}

\begin{document}
\title{Multi-product splitting and Runge-Kutta-Nystr\"om integrators}

\author{Siu A. Chin}

\affiliation{Department of Physics, Texas A\&M University,
College Station, TX 77843, USA}

\begin{abstract}
The splitting of $\e^{h(A+B)}$ into a single product 
of $\e^{h A}$ and $\e^{hB}$ results in symplectic integrators 
when $A$ and $B$ are classical Lie operators. 
However, at high orders, a single product splitting, with exponentially growing 
number of operators, is very difficult to derive. This work shows that, 
if the splitting is generalized to a sum of products, then a simple choice 
of the basis product reduces the problem to that of extrapolation, 
with analytically known coefficients and only quadratically growing 
number of operators. When a multi-product splitting is applied to 
classical Hamiltonian systems, 
the resulting algorithm is no longer symplectic but is of the 
Runge-Kutta-Nystr\"om (RKN) type. Multi-product splitting, 
in conjunction with a special force-reduction process,
explains why at orders $p=4$ and 6, RKN integrators only need
$p-1$ force evaluations.

\end{abstract}
\maketitle

\section {Introduction}

The approximation of 
${\rm e}^{\ep(A+B)}$ to any order in $\ep$ via
a single product decomposition
 \be
{\rm e}^{\ep( A+ B)}=
\prod_{i}
{\rm e}^{a_i\ep A}{\rm e}^{b_i\ep B}
\label{prod} 
\ee
is the basis of the splitting method\cite{mcl02,hairer02,lr04}
for solving diverse classical\cite{dragt,ruth83,neri87,cre89,fr90,suzu90,yos90,wis91,yos93},
quantum mechanical\cite{glasner,shen,serna,chinc01}, and stochastic 
evolution equations\cite{drozdov,for01}. The resulting algorithms are then
symplectic, unitary and norm-preserving, respectively. However, at high orders,
a single product splitting of the form (\ref{prod}) is very difficult to derive and 
the number of operators grows exponentially. This is easy to see: the error terms 
of the decomposition are governed by the fundamental Baker-Campbell-Hausdorff formula, 
\begin{equation}
{\rm e}^{\ep A}{\rm e}^{\ep B}=
{\rm e}^{\ep (A+B)+{1\over 2}\ep^2[A,B]+{1\over{12}}\ep^3([A,[A,B]]-[B,[A,B]])+\cdots}
\label{bch}
\end{equation}
where $[A,B]=AB-BA$ is the commutator bracket.
Starting with the first error term $[A,B]$, successive higher order error terms are 
generated by further commutators of $A$ and $B$, such as $[A,[A,B]]$, $[B,[A,B]]$,
 $[A,[A,[A,B]]]$ and $[B,[A,[A,B]]]$, etc..
This rapid doubling is
moderated somewhat by the Jacobi identity
and the special form of operators $A$ and $B$, but the growth of these
error terms is essentially exponential. 
This implies that high-order  
symplectic integrators must require many force-evaluations to eliminate these 
many error terms. For symplectic integrators, the minimum number of force-evaluations 
necessary at orders 4, 6, 8, 10 has been found empirically to be 3 (Ref.\cite{cre89,fr90,suzu90,yos90}), 
7 (Ref.\cite{yos90}), 15 (Ref.\cite{suzu93}) and 31 (Ref.\cite{kahan97}), respectively. 
The doubling trend is clearly visible and seems to be precisely
twice plus one. 

While symplectic integrators have many desirable conserving properties, such as
the exact conservation of Poincar\'e invariants, they are not immune from
the irreversible phase errors\cite{shita,gladman,chin2000,chin071} 
that directly affects the accuracy of  
trajectories. Since the accuracy of the trajectory is of paramount importance, 
regardless of how well some Porincar\'e invariants are preserved, 
it is less clear then, given the same number of force-evaluations, that 
a $2n$-order symplectic integrator should be preferred over a $2n+4$ or $2n+6$
order non-symplectic integrator, even for long term integrations.

This work shows that if one generalizes the decomposition of $\e^{h(A+B)}$ 
to a sum of products, then for a special choice of the basis product, the problem 
reduces to that of simple extrapolation and can be solved easily to 
any even order. When applied to classical Hamiltonian
systems, the resulting algorithms are no longer symplectic, but	
correspond to Runge-Kutta-Nystr\"om (RKN) type integrators. Blanes, Casas and
Ros\cite{bcr99} and Blanes and Casas\cite{bc05}
have previously shown that these extrapolated symplectic integrator enjoys better
conservation properties than conventional RKN algorithms.  
Interestingly, multi-product splitting 
not only reproduces Nystr\"om integrators at orders $p=4$ and 6 but
also explain why they only need $p-1$ force evalutions.
At even higher order, these extrapolated RKN algorithms are surprisingly 
competitive with symplectic and other RKN integrators found in the literature. 

\section {Multi-Product Expansion}

If ${\rm e}^{\ep(A+B)}$ were to be decomposed into a sum of products  
\be
{\rm e}^{\ep( A+ B)}=\sum_k c_k
\prod_{i}
{\rm e}^{a_{k,i}\ep A}{\rm e}^{b_{k,i}\ep B}
\label{mprod} 
\ee
then there is tremendous freedom in the choice of $\{c_k, a_{k,i},b_{k,i}\}$.
However, for the most general of applications, including
solving time-irreversible diffusion type equations\cite{drozdov,for01},  
one must keep the coefficients $\{a_{k,i},b_{k,i}\}$ positive.
This crucial and deliberate choice distinguishes this work from previous
extrapolations, which made no distinction between solving time-reversible and
time-irreversible problems.  
If $\{a_{k,i},b_{k,i}\}$ were to be positive, then by Sheng's theorem\cite{sheng89},
any such individual product can at most be second order and some 
coefficients ${c_k}$ must be negative. Thus without loss of generality,
one can choose the simplest basis product as either one of
the following symmetric second-order splittings  
\be
{\cal T}_2(\ep)=
{\rm e}^{\frac12 \ep B}
{\rm e}^{\ep A}
{\rm e}^{\frac12 \ep B}
\la{sec2a}
\ee
or
\be
{\cal T}_2(\ep)=
{\rm e}^{\frac12 \ep A}
{\rm e}^{\ep B}
{\rm e}^{\frac12 \ep A}.
\la{sec2b}
\ee
The choice of a symmetric product is important,
because it has only odd powers of $h$,
\be
{\cal T}_2(\ep)=\exp(\ep (A+B)+\ep^3 E_3+\ep^5 E_5+\cdots)
\la{secerr}
\ee
where $E_i$ are higher order error commutators of $A$ and $B$.
It then follows that the $k$th power of $\ct_2$ at step size $h/k$ is 
exactly given by
\be
{\cal T}_2^k\left(\frac\ep{k}\right) 
=\exp(\ep (A+B)+k^{-2}\ep^3 E_3+k^{-4}\ep^5 E_5+\cdots).
\la{thk}
\ee
Thus the set of powers ${\cal T}_2^{k}(\ep/k)$ 
has well known error structure and forms a suitable basis for 
expanding
${\rm e}^{\ep(A+B)}$. For example, any two terms with
$k=l$ and $k=m$ can approximate ${\rm e}^{\ep(A+B)}$ 
to fourth-order via 
\be
{\rm e}^{\ep(A+B)}
=c_l{\cal T}_2^{l}\left(\frac\ep{l}\right)
+c_m{\cal T}_2^{m}\left(\frac\ep{m}\right)+e_5(h^5E_5)
\la{fouro}
\ee
with obvious solutions
\be
c_l=\frac{l^2}{l^2-m^2}\qquad c_m=\frac{m^2}{m^2-l^2}
\ee
and error coefficient
\be
e_5=-\frac1{l^2m^2}.
\ee
The RHS of (\ref{fouro}) can also be written as
\be
{\cal T}_2^{l}\left(\frac\ep{l}\right)
+\frac{1}{(l/m)^2-1}\left[{\cal T}_2^{l}\left(\frac\ep{l}\right)
-{\cal T}_2^{m}\left(\frac\ep{m}\right)\right]
\ee
which coincide with the diagonal elements of the 
Richardson-Aitken-Neville extrapolation\cite{hair93} table.
Here, we do not do the extrapolation numerically, intead,
we give an analytical formula for extrapolating to any
even order. 
More generally, for a given set of $n$ distinct whole numbers $\{k_1, k_2, ...k_n\}$,
one can form a $2n$-order approximation of ${\rm e}^{\ep(A+B)}$ via 
\be
{\rm e}^{\ep(A+B)}
=\sum_{i=1}^n c_i{\cal T}_2^{k_i}\left(\frac\ep{k_i}\right)
+e_{2n+1}(h^{2n+1}E_{2n+1}).
\la{mulexp}
\ee
The expansion coefficients ${c_i}$ are determined by 
a specially simple Vandermonde equation, 
\be
\pmatrix{1     & 1        & 1      & ... & 1      \cr
 k_1^{-2}      & k_2^{-2}      & k_3^{-2} & ...      & k_n^{-2} \cr
 k_1^{-4}      & k_2^{-4}      & k_3^{-4} & ...      & k_n^{-4} \cr
      ...      & ...           & ...      & ...      & ...      \cr
 k_1^{-2(n-1)} & k_2^{-2(n-1)} & k_3^{-2(n-1)} & ... & k_n^{-2(n-1)} 
		 }
\pmatrix{c_1 \cr c_2\cr c_3\cr ...\cr c_n}
=\pmatrix{1\cr 0\cr 0\cr ...\cr 0 }
\la{solut}
\ee
with {\it closed form} solutions
\be
c_i=\prod_{j=1 (\ne i)}^n\frac{k_i^2}{k_i^2-k_j^2}
\la{coef}
\ee
and error coefficient,
\be
e_{2n+1}=(-1)^{n-1}\prod_{i=1}^n\frac1{k_i^2}.
\la{error}
\ee
The closed forms (\ref{coef}) and (\ref{error})
are the key results of this multi-product expansion. All error terms
of the same order, rather then be forced to zero individually, are extrapolated
to zero simultaneously. The Vandermonde equation (\ref{solut}) is a special
case of extrapolated symplectic algorithms considered by Blanes, Casas 
and Ros\cite{bcr99}. However, they did not obtain the exact solution given 
by (\ref{coef}). The proof that (\ref{coef}) is the exact solution is given 
in Appendix A. According to Blanes, Casas and Ros\cite{bcr99}, this class of 
extrapolated integrator is symplectic to order 2n+3, which at higher orders, 
can easily exceed one's machine precision. That is, at sufficiently high
orders, these extrapolated symplectic algorithms are, up to machine precision,
indistinguishable from truly symplectic algorithms.

Since a $2n$-order expansion is completely characterized by
a set of $n$ whole numbers $\{k_1, k_2, ... k_n\}$, the resulting 
algorithm will be referred to as the $\{k_1, k_2, ... k_n\}$-integrator. 
The error coefficient (\ref{error}) implies that, if the number of application 
of $\ct_2$ is fixed,
\be
\sum_{i=1}^n k_i=K
\la{fixf}
\ee
then the optimal algorithm with the least error is given by a set of 
$n$ distinct whole numbers $\{k_i\}$ closest to INT($K/n$) satisfying (\ref{fixf}). 
On the other hand, at order $2n$, the minimum number of $\ct_2$ required is given by the 
natural sequence $\{k_i\}=\{1, 2, 3 \,...\, n\}$, corresponding to	approximating
${\rm e}^{\ep(A+B)}$ 
via 
\be
{\rm e}^{\ep( A+ B)}=\sum_{k=1}^n c_k {\cal T}_2^k\left(\frac\ep{k}\right)
+e_{2n+1}(h^{2n+1}E_{2n+1})
\la{mcomp}
\ee
with $n(n+1)/2$ applications of $\ct_2$ and an error coefficient of
\be
e_{2n+1}=(-1)^{n-1}\frac1{(n!)^2}.
\ee
The simplicity of these results make multi-product splitting 
algorithms extremely easy to analyze and implement. 

From (\ref{coef}), the minimum-$\ct_2$ expansions of orders
four to ten are given by 
\be
{\cal T}_4(\ep)=-\frac13{\cal T}_2(\ep)
+\frac43{\cal T}_2^2\left(\frac\ep{2}\right)
\la{four}
\ee
\be
{\cal T}_6(\ep)=\frac1{24} {\cal T}_2(\ep)
-\frac{16}{15}{\cal T}_2^2\left(\frac\ep{2}\right)
+\frac{81}{40}{\cal T}_2^3\left(\frac\ep{3}\right)
\la{six}
\ee
\be
{\cal T}_8(\ep)=-\frac1{360} {\cal T}_2(\ep)
+\frac{16}{45}{\cal T}_2^2\left(\frac\ep{2}\right)
-\frac{729}{280}{\cal T}_2^3\left(\frac\ep{3}\right)
+\frac{1024}{315}{\cal T}_2^4\left(\frac\ep{4}\right)
\la{eight}
\ee
\be
{\cal T}_{10}(\ep)=\frac1{8640} {\cal T}_2(\ep)
-\frac{64}{945}{\cal T}_2^2\left(\frac\ep{2}\right)
+\frac{6561}{4480}{\cal T}_2^3\left(\frac\ep{3}\right)
-\frac{16384}{2835}{\cal T}_2^4\left(\frac\ep{4}\right)
+\frac{390625}{72576}{\cal T}_2^5\left(\frac\ep{5}\right).
\la{ten}
\ee
In contrast to a single product decomposition, whose
coefficients $\{a_i,b_i\}$ are generally irrational and 
must be determined order-by-order numerically with limited precision, 
the multi-product expansion has only rational coefficients and
are given analytically by (\ref{coef}) for all even orders.

The operator extrapolation here is at once more general and simpler than
extrapolating solutions of differential equations. The general
expansion (\ref{mcomp}) can be applied to any evolution operator 
$\e^{h(A+B)}$ in terms of ${\cal T}_2$, with known second order
error structure (\ref{secerr}). There is no need to devise and 
prove\cite{gragg65}, a second-order solution in advance.
The low order operator extrapolation (\ref{four}) has been 
used previously\cite{kevin95,schatz02} in different contexts. 
Here, we provide a systematic expansion to any even order
and point out the connection between symplectic and RKN integrators.

\section {Deriving RKN integrators}

As emphasized in the last section, the multi-product
expansion (\ref{mcomp}) is an extrapolated operator
approximation to the evolution operator $\e^{h(A+B)}$
and can be applied to many different types of equations.
Its applications to to quantum\cite{chin093} and stochstic\cite{chin094} 
evolution equations have already proven to be highly successful.
Here, we will concentrate on its use in solving classical
dynamic	problems.

If ${\rm e}^{\ep(A+B)}$ is decomposed into a 
sum of products as in (\ref{mcomp}), then the resulting algorithm will no 
longer be symplectic. This is because for more than one product, cross 
terms will spoil the usual proof of symplecticity 
(See Ref.\cite{dragt}, Theorem 1.). As we will show, multi-product
splitting now produces RKN type integrators.

For solving Hamilton's equation in the form,
\be
\frac{d \bq}{dt}=\bv \qquad \frac{d\bv}{dt}=\bac(\bq)
\ee
the operators are
\be
A=\sum_i v_i\frac{d}{dq_i}  \qquad B=\sum_ia_i(\bq)\frac{d}{dv_i} 
\ee
and where we have abbreviated $\bv=\bp/m$ and $\bac(\bq)=\bff(\bq)/m$.
Thus (\ref{sec2a}) corresponds to the velocity-form of the Verlet algorithm (VV):
\ba
\bv_1&=&\bv_0+\frac{h}2\,  \bac(\bq_0)\cr
\bq_1&=&\bq_0+h\, \bv_1	\cr
\bv_2&=&\bv_1+\frac{h}2\, \bac(\bq_1)
\la{vv}
\ea
and (\ref{sec2b}) corresponds to the
position-form of the Verlet algorithm (PV):
\ba
\bq_1&=&\bq_0+ \frac{h}2\,\bv_0	\cr
\bv_1&=&\bv_0+h\, \bac(\bq_1)\cr
\bq_2&=&\bq_1+ \frac{h}2\,\bv_1.
\la{pv}
\ea
The last numbered variables are the updated variables. Both VV and PV
are second order symplectic integrators. The operator
$\ct_2^k(h/k)$ simply iterating either algorithm $k$ times 
at step size $h/k$ and is therefore also symplectic. 
The PV algorithm uses only one force
evaluation per update. If $\ct_2$ is taken to be PV, then
the extrapolated integrators of order $2n$, as examplified by
(\ref{four})-(\ref{ten}), only require $n(n+1)/2$ force-evaluations.
The minimal integrator of orders 4, 6, 8 and 10 requires only 3, 6, 10, and 
15 force evaluations, respectively. 

If $\ct_2$ is taken to the VV algorithm, since all products use
the same starting force, each $\ct_2^k(h/k)$ for $k>2$ also 
uses only $k$ force-evaluations. However, $\ct_2(h)$ must then evaluate
the initial force for the rest of the product and its
own final force. Thus $\ct_2(h)$ requires two force-evaluations, giving
the final force count as $n(n+1)/2+1$. 

There are other possible extrapolation methods, such as
\be
{\cal T}_6(\ep)=-\frac1{15}{\cal T}_4(\ep)
+\frac{16}{15}{\cal T}_4^2\left(\frac\ep{2}\right)
\la{sixp}
\ee
which is similar to the triplet concatenation for symplectic 
integrators\cite{cre89,fr90,suzu90,yos90}.
However, such an iterative extrapolation, which triples the number of 
force-evaluations in going from order $2n$ to $2n+2$, is not competitive 
with (\ref{mcomp})'s linear increase of only $n+1$ additional force-evaluations.

To see how the expansion (\ref{mcomp}) results in  
RKN integrators, let's takes 
$\ct_2$ to be VV, then the extrapolation (\ref{four}) produces the 
following fourth-order, $\{1,2\}$-integrator
\ba
&&\bq=\bq_0+\ep\,\bv_0+\frac{\ep^2}6\left(\bac_0+2\bac_{1/2}\right)\nn\\
&&\bv=\bv_0+\frac{\ep}6\left(\bac_0+4\bac_{1/2}
+2\bac_{2/2}-\bac_{1/1}\right)
\la{rknn}
\ea
where we have systematically denoted the force evaluation point of $\ct^k_2(h/k)$
at the intermediate step $(i/k)h$ as $\bq_{i/k}$ and the resulting force
as $\bac_{i/k}=\bac(\bq_{i/k})$, with $\bac_0=\bac(\bq_0)$. 
The final force evaluation point of $\ct_2(h)$, the midpoint of $\ct_2^2(h/2)$, 
and the final position of 
$\ct_2^2(h/2)$ are given by
\ba
&&\bq_{1/1}=\bq_0+\ep\bv_0+\frac{\ep^2}{2}\bac_0\nn\\
&&\bq_{1/2}=\bq_0+\frac\ep{2}\bv_0+\frac{\ep^2}{8}\bac_0\nn\\
&&\bq_{2/2}=\bq_0+\ep\bv_0+\frac{\ep^2}{4}\left(\bac_0+\bac_{1/2}\right).
\la{rknq}
\ea
The resulting integrator (\ref{rknn}) appears to require four force-evaluations,
in accordance with the formula $n(n+1)/2+1$. However, a key observation here is   
that force subtractions
at the {\it same time step} can be combined into a single force evaluation. 
Let
\be
\delta\bq_2= \bq_{2/2}-\bq_{1/1}=\frac{h^2}4 (\bac_{1/2}-\bac_0)\approx O(h^3)
\la{dq2}
\ee
then the subtraction of two forces gives,
\ba
2\bac_{2/2}-\bac_{1/1}
&=&2\bac(\bq_{1/1}+\delta\bq_2)-\bac(\bq_{1/1})\la{coll}\\
&=&\bac( \bq_{1/1}+2\delta\bq_2 )+O(h^6)\nn\\
&=&\bac\left(\bq_0+\ep\bv_0+\frac{1}{2}\ep^2\bac_{1/2}\right)+O(h^6).
\la{tqo}
\ea
The total number of force evaluations is reduced from four to three and 
(\ref{rknn}) reproduces Nystr\"om's original fourth-order integrator\cite{nys25,bat99}.
Thus we have shown analytically that, the subtraction of
two {\it symplectic} integrators reproduces Nystr\"om's original integrator.
This connection between extrapolated symplectic integrators and Nystr\"om's integrators
has not been appreciated previously.     
 
If $\ct_2(h)$ is taken to be PV, the corresponding $\{1,2\}$-integrator
with three force-evaluations is given by
\ba
&&\bx=\bx_0+\ep\,\bv_0+\frac{\ep^2}6\left(3\bac_{1/4}-\bac_{1/2}+\bac_{3/4}\right)\nn\\
&&\bv=\bv_0+\frac{\ep}3\left(2\bac_{1/4}-\bac_{1/2}
+2\bac_{3/4}\right).
\la{rknpv}
\ea
Since PV and VV based algorithms have different force-evaluation points, 
to avoid confusion, we will 
denote the positions of PV-based integrator as $\bx_{i/k}$. 
Here, the force-evaluation points are
\ba
&&\bx_{1/4}=\bx_0+\frac{\ep}4 \bv_0\nn\\
&&\bx_{1/2}=\bx_0+\frac\ep{2}\bv_0\nn\\
&&\bx_{3/4}=\bx_0+\frac34 \ep\bv_0+\frac{\ep^2}{4}\bac_{1/4}.
\la{rknpvx}
\ea
All such PV-based integrators are non-FSAL (First Step As Last) RKN integrators 
in that the force is never evaluated initially. The explicit
force subtraction will be 
more susceptible to round-off errors. However, this fourth-order
integrator is superior to the Nystr\"om	integrator when solving
Kepler's orbit and concerns for round-off errors are 
moderated by the wide use of quadruple and multi-precision packages.	 

To compare the long time integration error of numerical algorithms,
it is useful to solve Kepler's orbit in two dimension,
\be
a_x=-\frac{q_x}{q^3}\quad 
a_y=-\frac{q_y}{q^3}\quad q=\sqrt{q_x^2+q_y^2}
\ee
with initial conditions
\be
v_x=0\quad v_y=\sqrt{\frac{1-e}{1+e}}\quad q_x=1+e\quad q_y=0
\ee
where $e$ is the eccentricity of the orbit.
The resulting energy and period are respectively $-1/2$ and $2\pi$.
This set of initial conditions is chosen (instead of the usual one in
the literature) because it is non-singular as $e\rightarrow 1$ and the
orbit's semi-major axis is always along the x-axis for all values of $e$.
Both are useful for computing the precession error of 
Kepler's orbit at high eccentricity.  

To gauge the accuracy of the trajectory intrinsically, without
any external comparison, we monitor the irreversible 
phase error\cite{shita,gladman} by computing the rotation angle per period, 
$\Delta\theta$, of 
the Laplace-Runge-Lenz (LRL) vector\cite{chin2000,chin071}:
\be
\A=\bv\times{\bf L}-\frac{\bq}q
\ee
with ${\bf L}=\bq\times\bv$. 
If the orbit is exact, then the LRL vector would remain 
constant pointing along the semi-major axis of the orbit. If 
$\Delta\theta\ne 0$, then the orbit will have precessed an angle
of $\Delta\theta$ after each period. (This is the advance of the
perihelion of solar planets.) Thus $\Delta\theta$ is a much more
direct measure of the trajectory's accuracy than the energy error.
This precession error grows linearly with time as $m\Delta\theta$, where $m$ is the
number of periods, even for symplectic integrators\cite{shita,gladman}. The 
corresponding precession error coefficient is extracted by dividing $\Delta\theta$ by
$h^4$ using smaller and smaller $h$ until 
\be
\lim_{h\rightarrow 0}\, \frac{\Delta\theta}{h^{4}}=e_P
\la{errorpre}
\ee
is independent of $h$. The coefficient $e_P$ as a function
of the eccentricity is therefore a unique error signature of each 
fourth-order integrator independent of the step size. 

\begin{figure}
	\centerline{\includegraphics[width=0.8\linewidth]{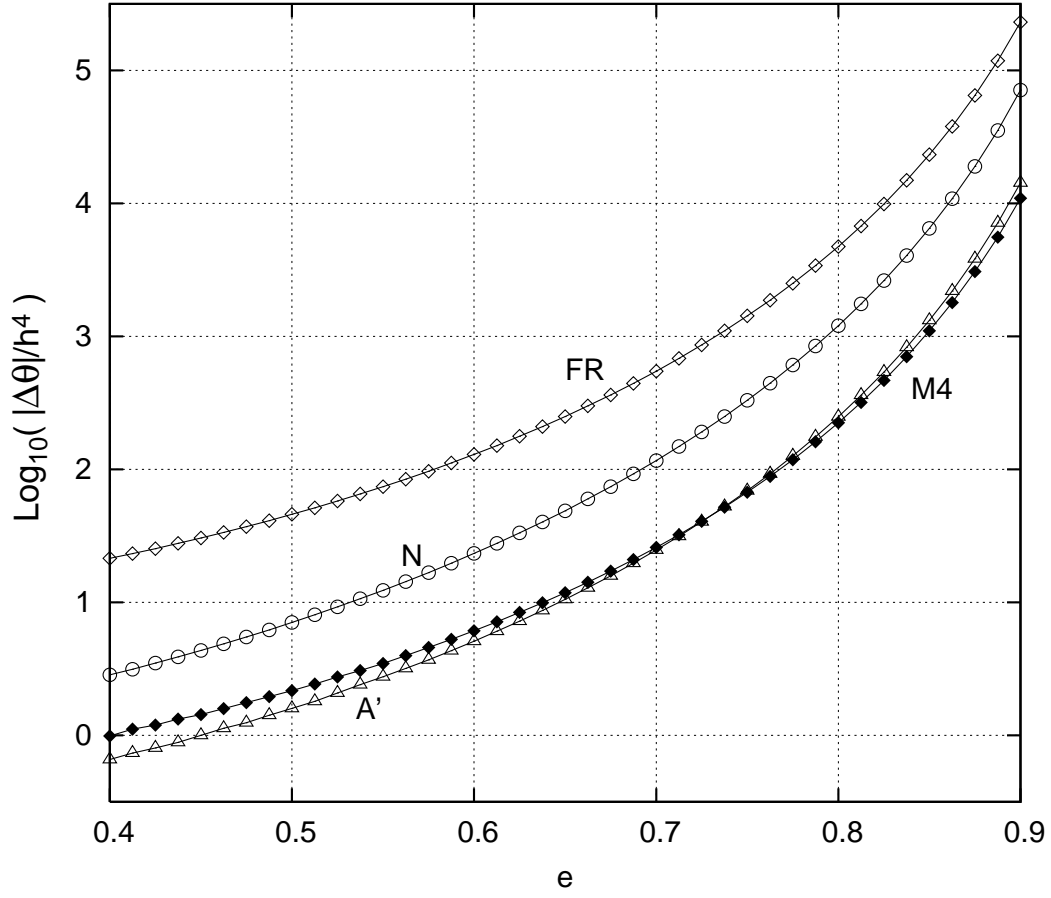}}
\caption{
The orbital precession error coefficient as a function of the eccentricity
$e$ of the Kepler orbit for four fourth-order integrators which use 
only three force evaluations per update. FR is the
Forest-Ruth symplectic integrator, N is the original 
Nystr\"om integrator,
A$^\prime$ is a forward integrator and M4 is the multi-product splitting
integrator (\ref{rknpv}). 
\label{fourecc}}
\end{figure}

In Fig.\ref{fourecc}, we compare the precession error coefficient of four distinct 
fourth-order integrators with only three force evaluations. The symplectic Forest-Ruth (FR)
integrator\cite{fr90}, the forward integrator A$^\prime$\cite{ome06,chin072}, the Nystr\"om integrator (N) given
by (\ref{rknn}), (\ref{tqo}) and the minimal fourth-ordr $\{1,2\}$-integrator M4, given by
(\ref{rknpv}). The forward integrator A$^\prime$ has only 
positive splitting coefficients $a_i$ and $b_i$. The convergence of (\ref{errorpre})
is seen near $h=2\pi/3000$; the final value used is $h=2\pi/5000$. 
The coefficient $e_P$ grows steeply with $e$, by four orders of magnitude
from $e=0.4$ to $e=0.9$, but all integrators showed similar dependence on $e$.
It is well known that the precession error of the symplectic FR integrator 
is greater than that of RKN algorithms\cite{shita,gladman} such as N, but it
was not known previously that M4 is so much better than FR and N. 
The error coefficient $e_P$ at $e=0.9$ for FR, N, A$^\prime$ and M4
are -23.1$\times10^4$, 7.1$\times10^4$, -1.4$\times10^4$ and -1.1$\times10^4$, respectively.

\section {Sixth Order RKN integrators}

If $\ct_2(h)$ is taken to be VV, the difference between the
two force evaluation points at the same intermediate time step, such
as (\ref{dq2}), is always of order $O(h^3)$ and the resulting force subtration, 
such as (\ref{tqo}), is of order $O(h^6)$. This immediately suggests that
the three final force evaluations
in the $\{1,2,3\}$-integrator can be collapsed into one, yielding a
sixth order algorithm with only five force-evaluations. This is indeed the 
case. The $\{1,2,3\}$-integrator, according to (\ref{six}), is given by
\be
\bq=\bq_0+\ep\,\bv_0
+\frac{\ep^2}{120}\bigl(
11\bac_0
+54\bac_{1/3}
-32\bac_{1/2}
+27\bac_{2/3}
\bigr),\la{al6123q}
\ee
\ba
\bv=\bv_0+\frac{\ep}{240}\bigl(
 22\bac_0
+162\bac_{1/3}
&-&128\bac_{1/2}
+162\bac_{2/3}\nn\\
&+& 81\bac_{3/3}
-64\bac_{2/2}
+5\bac_{1/1}\bigr),
\la{al6123p}
\ea
with three additional force evaluations at:
\ba
&&\bq_{1/3}=\bq_0+\frac\ep{3}\bv_0+\frac{\ep^2}{18}\bac_0,\\
&&\bq_{2/3}=\bq_0+\frac2{3}\ep\bv_0+\frac{\ep^2}{9}\left(\bac_0+\bac_{1/3}\right),\\
&&\bq_{3/3}=\bq_0+\ep\bv_0+\frac{\ep^2}{18}\left(3\bac_0
+4\bac_{1/3}+2\bac_{2/3}\right).\la{q33q}
\ea
Let
\be
\delta\bq_3=\bq_{3/3}-\bq_{1/1}=\frac1{9} h^2(2\bac_{1/3}+\bac_{2/3}-3\bac_0)\approx O(h^3),
\la{dq3}
\ee
then the three final-step force-evaluations can be collapsed into one,
\be
 81\bac_{3/3}-64\bac_{2/2}+5\bac_{1/1}
 =22\bac(\widetilde\bq_1)+O(h^6),
 \la{fcol}
\ee
with
\ba
\widetilde\bq_1&=&\bq_{1/1}+\frac1{22}(81\delta\bq_3-64\delta\bq_2)\nn\\
&=&\bq_0+\ep\bv_0+\frac{\ep^2}{22}\left(
18\bac_{1/3}-16\bac_{1/2}+9\bac_{2/3}\right).
\la{qst}
\ea
The force consolidation (\ref{fcol}) now renders (\ref{al6123p}) symmetric, 
\be
\bv=\bv_0+\frac{\ep}{240}\left(
 22\bac_0
+162\bac_{1/3}
-128\bac_{1/2}
+162\bac_{2/3}
+22\widetilde\bac_1
\right).
\la{al6123pp}
\ee
This sixth order integrator seems not to be known prior to this work. 
The well known
sixth order integrator with five force evaluations due to Albrecht\cite{alb55},
can be derived from an alternative expansion,
\be
{\cal T}_6(\ep)=\frac1{45} {\cal T}_2(\ep)
-\frac{4}{9}{\cal T}_2^2\left(\frac\ep{2}\right)
+\frac{64}{45}{\cal T}_2^4\left(\frac\ep{4}\right),
\la{sixalb}
\ee
corresponding to the integrator $\{1,2,4\}$,
\be
\bq=\bq_0+\ep\,\bv_0
+\frac{\ep^2}{90}\bigl(
7\bac_0
+24\bac_{1/4}
+16\bac_{2/4}
-10\bac_{1/2}
+8\bac_{3/4}\bigr),
\la{al6q}\ee
\ba
\bv=\bv_0+\frac{\ep}{90}\bigl(7\bac_0
+32\bac_{1/4}
+32\bac_{2/4}
&-& 20\bac_{1/2}
+32\bac_{3/4}\nn\\
&+&16\bac_{4/4}
-10\bac_{2/2}
  +\bac_{1/1}\bigr),
\la{al6p}
\ea
with four additional force evaluations at
\ba
&&\bq_{1/4}=\bq_0+\frac\ep{4}\bv_0+\frac{\ep^2}{32}\bac_0,\\
&&\bq_{2/4}=\bq_0+\frac\ep{2}\bv_0+\frac{\ep^2}{16}\left(\bac_0+\bac_{1/4}\right),\\
&&\bq_{3/4}=\bq_0+\frac34 \ep\bv_0+\frac{\ep^2}{32}\left(3\bac_0
+4\bac_{1/4}+2\bac_{2/4}\right),\la{q3q}\\
&&\bq_{4/4}=\bq_0+\ep\bv_0+\frac{\ep^2}{16}\left(2\bac_0
+3\bac_{1/4}+2\bac_{2/4}+\bac_{3/4}\right).
\la{rkn4q}
\ea
This algorithm nominally requires eight force evaluations, but
forces can now be collapsed at $h/2$ and at $h$. At $h/2$, one has
\be
16\bac_{2/4}-10\bac_{1/2}
=6\bac(\bq^*_{1/2})+O(h^6),
\la{chalf}
\ee
where the new half time-step evaluation point is
\be
\bq^*_{1/2}
=\bq_0+\frac\ep{2}\bv_0+\frac{\ep^2}{24}\left(4\bac_{1/4}-\bac_0\right).
\ee
Similarly, the three force evaluations at $h$
can be collapsed into one,  
\be
16\bac_{4/4}
-10\bac_{2/2}
+\bac_{1/1}=7\bac(\bq^\dagger_1)+O(h^6),
\la{cone}
\ee
where
\be
\bq^\dagger_1=\bq_0+\ep\bv_0+\frac{\ep^2}{14}\left(
6\bac_{1/4}+4\bac_{2/4}-5\bac_{1/2}+2\bac_{3/4}\right).
\la{qsone}
\ee
The number of force-evaluations
can be reduced from eight to five if both $\bq_{2/4}$ 
and $\bq_{1/2}$ can be replaced everywhere by
$\bq^*_{1/2}$. However, replacing $\bac_{2/4}$ in $\bq_{3/4}$ by $\bac^*_{1/2}$,
would yield
\ba
\bq^*_{3/4}&=&\bq_0+\frac34 \ep\bv_0+\frac{\ep^2}{32}\left(3\bac_0
+4\bac_{1/4}+2\bac^*_{1/2}\right)\nn\\
&=&\bq_{3/4}+\frac{\ep^2}{16}\left(\bac^*_{1/2}-\bac_{2/4}\right)\nn\\
&=&\bq_{3/4}+O(h^5),
\la{star3q}
\ea
and the resulting $32\bac^*_{3/4}$ term	in (\ref{al6p}) would produce an error of 
$O(h^6)$, spoiling the algorithm. 
Remarkably, replacing $\bq_{2/4}$ 
and $\bq_{1/2}$ everywhere by $\bq^*_{1/2}$ in (\ref{qsone}) defines
\ba
\bq^*_1&=&\bq_0+\ep\bv_0+\frac{\ep^2}{14}\left(
6\bac_{1/4}-\bac^{*}_{1/2}+2\bac^{*}_{3/4}\right)\nn\\
&=&\bq^\dagger_1-\frac{\ep^2}{14}\left(
\bac^{*}_{1/2}+4\bac_{2/4}-5\bac_{1/2}
\right), \nn\\
&=&\bq^\dagger_1+O(h^5),
\la{qsonep}
\ea
and the resulting $O(h^6)$ error of $7\bac^*_1$ 
in (\ref{al6p}) exactly cancels that of $32\bac^*_{3/4}$! Thus,
one recovers Albrecht's integrator:
\ba
&&\bq=\bq_0+\ep\,\bv_0+\frac{\ep^2}{90}\left(7\bac_0+24\bac_{1/4}
+6\bac^*_{1/2}+8\bac^*_{3/4}\right)+O(h^7)
\la{al6qp}\\
&&\bv=\bv_0+\frac{\ep}{90}\left(7\bac_0
+32\bac_{1/4}
+12\bac^{*}_{1/2}
+32\bac^{*}_{3/4}
+7\bac^{*}_{1}
\right)+O(h^7).
\la{al6pp}
\ea
Within the context of multi-product expansion,
these are the only two sixth-order integrators possible with only
five force evaluations. Again, we have shown analytically that the extrapolation
of three second-order symplectic integrators produces sixth-order RKN integrators.

\begin{figure}
	\centerline{\includegraphics[width=0.8\linewidth]{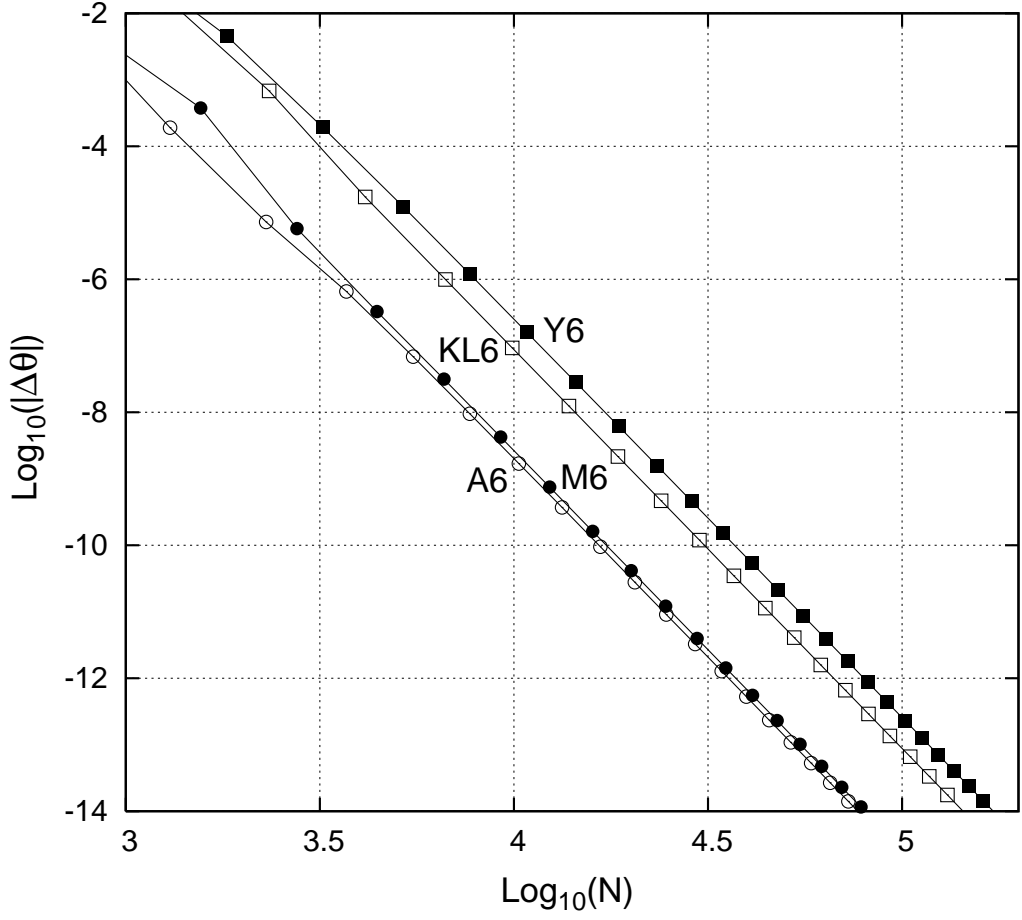}}
\caption{
The orbital precession error of sixth-order 
integrators as a function of
of the number of force evaluation $N$. Y6, KL6 and A6 are integrator
of Yoshida, Kahan-Li and Albrecht respectively. M6 is the minimal extrapolated
integrator corresponding to (\ref{six}).  
\label{sixecc}}
\end{figure}

To compare integrators of the same order but with different number of force 
evaluations, we compute the precession error $\Delta\theta$ at $e=0.9$
as a function of $N$, the number of force-evaluations used.
The results for four sixth-order algorithms are shown in Fig.\ref{sixecc}. 
Y6 is Yoshida's symplectic algorithm\cite{yos90} with the minimum 7 force evaluations.
KL6 is a much improved version by Kahan and Li\cite{kahan97} with 9 force evaluations. 
A6 and M6 are Albrecht's and the PV-based, minimal \{1,2,3\}-integrator  
with 5 and 6 
force-evaluations respectively. While KL6's error is only half of Y6's 
at $10^5$ force-evaluations, their respective error are
nearly 50 and 100 times larger than those of RKN integrators A6 and M6.
This calculation was carried out in quadruple-precision using the widely available
free software by Miller\cite{miller}.

\section {Higher order comparisons}

Since the force reduction process considered is  
only of $O(h^6)$, this process will no longer be effective  
for integrators beyond sixth-order. Thus beyond sixth-order, some other
force reduction processes must be found if 
the number of force-evaluation is to be less than $n(n+1)/2$. 
At higher orders, since the expansion coefficients (\ref{coef}) 
are known explicitly, it is trivial to write a subroutine for (\ref{pv}) 
and simply call it according the multi-product expansion (\ref{mcomp}). 
Alternatively, one can also construct them from the Aitken-Neville table\cite{hair93},
as done conventionally. It is of interest to compare these extrapolated
RKN integrators at high orders with existing symplectic and specially 
derived RKN integrators found in the literature.

\begin{figure}
	\centerline{\includegraphics[width=0.8\linewidth]{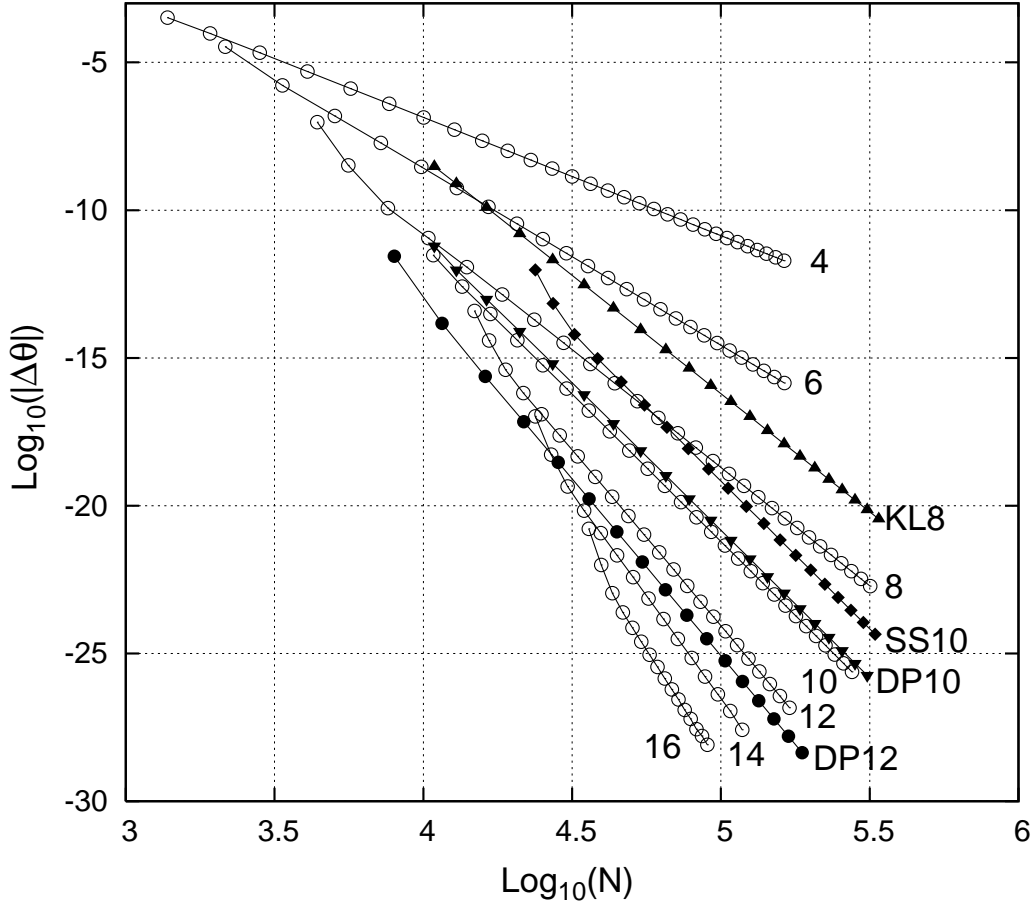}}
\caption{
The orbital precession error as a function of the number of force
evaluation $N$. KL8 and SS10 are Kahan-Li's eighth-order and
Sofroniou-Spalletta's tenth-order symplectic integrators respectively.
DP10 and DP12 are the 10th and 12th order RKN pair of Dormand and Prince.
The numbers 4, 6, 8, etc., denote the $2n$-order
integrators corresponding to the minimal multi-product expansion (\ref{mcomp}). 
\label{all}}
\end{figure}

In Fig.\ref{all}, we compare some well known higher order integrators
with the PV-based, minimal extrapolated integrators (\ref{mcomp}).
KL8 and SS10 are eighth and tenth order
symplectic integrators by Kahan and Li\cite{kahan97} and 
Sofroniou and Spalletta\cite{sof} with 17 and 35 force evaluations
respectively. Both are recommended by Ref.\cite{hairer02}. 
Because of their large number of force-evaluations, they are not 
competitive\cite{calvo93} in a precision-effort comparison as shown in Fig.\ref{all}.
At $N=10^5$, the error of KL8 is more than 300 times that of the
PV-based $\{1,2,3,4\}$-integrator denoted simply as 8 rather than M8. 
Similarly, the error of SS10 at $N=10^5$ is about 100 times that of
the minimal $\{1,2,3,4,5\}$-integrator denoted as 10.
DP10 and DP12 are the RKN integrator-pair by Dormand, El-Mikkawy and Prince\cite{dorm87}
with 17 force evaluations.
The coefficients for this integrators-pair were taken from 
Brankin {\it et al.}\cite{brank89} and converted to quadruple precision. 
DP10 is comparable to 10, which uses 15 force evaluations.
DP12's error at $N=10^5$ is 10 times smaller than 12. However, DP12's
superior performance can be quickly matched by going to integrator 14, whose error
is then 100 times smaller than DP12 at $N=10^5$. 
 
Fig.\ref{all} implies that the most efficient integrator for attaining a given
level of precision is never a fixed order integrator. As
one demands greater and greater precision, one must increase the order of the
integrator accordingly. In Fig.\ref{all}, the error curve of the $2n+2$ 
integrator is plotted only when it is below that of the
$2n$ integrator. The enveloping error curve is steeper than
any fixed order integrator. This maximum efficiency can be achieved 
only when one has the mean of producing an arbitrary high-order integrator 
at will, as it is done here. The highest order considered
in Fig.\ref{all} is 16 because quadruple precision is inadequate for precision below
$10^{-30}$. 

\section {Summary and conclusions}

In this work, we have generalized the single product decomposition of
the evolution operator $\e^{h(A+B)}$ to that of a sum of products.
By using powers of ${\cal T}_2$ as basis products, the multi-product expansion 
reduces to that of a simple extrapolation, with analytically known solutions.
Because the extrapolation is formulated on the level of operators, it
can be applied with great generality in solving diverse evolution
equations not only in classical mechanics but also in quantum and 
stochastic dynamics.
 
A multi-product expansion loses some desirable properties, such as not
being symplectic, unitary, norm preserving, etc., as compared to a single
product decomposition. However, at high orders, it gains in economy,
needing only a quadratic, rather than an exponential number of
operators. This is a vital consideration when high precision and
high order results are needed.

In applying to classical Hamiltonian dynamics, this work showed that
extrapolating symplectic integrators produces RKN integrators.
The process of force consolidation provided a simple explanation for
why $p-1$ force evaluations are sufficient for RKN integrators at orders 
$p=4$ and 6. An order-barrier then naturally arises when forces can longer be 
consolidated at higher $p$. Moreover, this work showed that, 
there exists a sequence of easily implemented, $2n$-order extrapolated RKN integrators, 
with only $n(n+1)/2$ force evaluations. Thus any other RKN integrator with fewer 
force evaluations is special and must embody some unique force consolidation 
schemes. These extrapolated RKN integrators are surprisingly competitive
with existing symplectic and specially devised RKN algorithms in
solving the Kepler problem. Since this series of extrapolated RKN integrators 
can be easily produced by anyone, they can serve as a common branchmark by which 
more sophisticated high order integrators can be judged.
For example, of all the integrators compared in this work, only
DP12 has outperformed its corresponding extrapolated RKN integrator
by having fewer force evaluations (17) than of $n(n+1)/2$ (=21).
Thus DP12 must possess some highly nontrivial force consolidation schemes.


\appendix
\section {Proving the solution}

To solve the Vandermonde equation (\ref{solut}) for $\{c_i\}$, 
let $x_i=k_i^{-2}$ so that (\ref{solut}) reads conventionally
\be
\pmatrix{1     & 1        & 1      & ... & 1      \cr
 x_1        & x_2      & x_3 & ...      & x_n \cr
 x_1^2      & x_2^2      & x_3^2 & ...      & x_n^2 \cr
      ...      & ...           & ...      & ...      & ...      \cr
 x_1^{(n-1)} & x_2^{(n-1)} & x_3^{(n-1)} & ... & x_n^{(n-1)} 
		 }
\pmatrix{c_1 \cr c_2\cr c_3\cr ...\cr c_n}
=\pmatrix{1\cr 0\cr 0\cr ...\cr 0 }
\la{solutx}
\ee
Consider the usual Lagrange interpolation at $n$ points $\{x_1, x_2, \dots, x_n\}$
with values $\{y_1, y_2, \dots, y_n\}$. The interpolating $n-1$ degree polynomial is given
by
\be
f(x)=\sum_{i=1}^n y_i L_i(x),
\la{lpoly}
\ee
where $L_i(x)$ are the Lagrange polynomials given by
\be
L_i(x)=\prod_{j=1 (\ne i)}^n \frac{(x-x_j)}{(x_i-x_j)}.
\ee
Since by construction
\be
L_i(x_k)=\delta_{ik}
\ee
the interpolation polynomial (\ref{lpoly}) correctly gives,
\be
f(x_k)=\sum_{i=1}^n y_i L_i(x_k)=\sum_{i=1}^n y_i \delta_{ik}=y_k.
\ee
Now let $y_i=x_i^m$ for $0\le m\le n-1$, 
then the interpolating polynomial
\be
f(x)=\sum_{i=1}^n x_i^m L_i(x)
\ee
and the function
\be
g(x)=x^m
\ee
both interpolate the same set of points. Since interpolating polynomials
of the same order
are unique, we must have $f(x)=g(x)$ and hence
\be
\sum_{i=1}^n x^m_iL_i(x)=x^m.
\ee
Writing this out in matrix form yields
\be
\pmatrix{1     & 1        & 1      & ... & 1      \cr
 x_1        & x_2      & x_3 & ...      & x_n \cr
 x_1^2      & x_2^2      & x_3^2 & ...      & x_n^2 \cr
      ...      & ...           & ...      & ...      & ...      \cr
 x_1^{(n-1)} & x_2^{(n-1)} & x_3^{(n-1)} & ... & x_n^{(n-1)} 
		 }
\pmatrix{L_1(x) \cr L_2(x)\cr L_3(x)\cr ...\cr L_n(x)}
=\pmatrix{1\cr x\cr x^2\cr ...\cr x^{n-1} }
\la{solutl}
\ee
Compare this to (\ref{solutx}), we immediately see that the solution is
\ba
c_i&=&L_i(0)\nn\\
&=&\prod_{j=1 (\ne i)}^n \frac{x_j}{x_j-x_i}\nn\\
&=&\prod_{j=1 (\ne i)}^n \frac{k_i^{2}}{k_i^{2}-k_j^{2}}
\ea
This proof is a simple application of the more general result of
inverting the Vandermonde matrix\cite{mem03}.


\newpage
\newpage
\newpage

\end{document}